\documentclass{amsart}
\usepackage{tikz}
\usepackage{hyperref}   
\usepackage{amssymb,xcolor}
\usepackage{stmaryrd} 

\usepackage{mathtools}  
\mathtoolsset{showonlyrefs}

\theoremstyle{plain}
\newtheorem{theorem}{Theorem}[section]
\theoremstyle{plain}

\theoremstyle{plain}
\newtheorem{definition}[theorem]{Definition}
\theoremstyle{plain}
\newtheorem{lemma}[theorem]{Lemma}
\theoremstyle{remark}
\newtheorem{remark}[theorem]{Remark}
\theoremstyle{plain}
\newtheorem{proposition}[theorem]{Proposition}
\numberwithin{equation}{section}
\theoremstyle{plain}

\begin{document}

\title{A tilted spacetime positive mass theorem}

\author{Xiaoxiang Chai}
\address{Department of Mathematics, POSTECH, Pohang, Gyeongbuk, South Korea}
\email{xxchai@kias.re.kr, xxchai@postech.ac.kr}

\begin{abstract}
  We show a spacetime positive mass theorem for asymptotically flat initial
  data sets with a noncompact boundary. We develop a mass type invariant and a
  boundary dominant energy condition. Our proof is based on spinors.
\end{abstract}



{\maketitle}

\section{Introduction}

The first positive mass theorem was proved by Schoen and Yau in their seminal
work {\cite{schoen-proof-1979}} using a minimal surface technique. It says
that if a complete manifold which is asymptotically flat and with nonnegative
scalar curvature, an quantity called the ADM mass defined at infinity is
nonnegative. The ADM mass is a characterization of scalar curvature at
infinity. There are various works on the positive mass theorem:
{\cite{witten-new-1981}}, {\cite{eichmair-spacetime-2016}},
{\cite{andersson-rigidity-2008}}, {\cite{wang-mass-2001}},
{\cite{chrusciel-mass-2003}}, {\cite{sakovich-jang-2021}}. Here the list is by
no means exhaustive.

The study of the positive mass type theorems of the asymptotically flat
manifold with a noncompact boundary were started in
{\cite{almaraz-positive-2016}}. As a result, the effect of the mean curvature
was included to the infinity and a boundary term was added to the ADM mass.
See {\cite{almaraz-mass-2020}}, {\cite{almaraz-spacetime-2019}},
{\cite{chai-positive-2018}}, {\cite{chai-asymptotically-2021-arxiv}} and
{\cite{almaraz-rigidity-2022-arxiv}} for some developments to the spacetime
and hyperbolic settings.

We revisit the asymptotically flat initial data set with a noncompact boundary
in this paper, we introduce the boundary dominant condition \eqref{tilt dec}
and we prove two spacetime positive mass theorems (Theorems \ref{tilt pmt},
\ref{tangential pmt}). We use the spinorial argument of Witten
{\cite{witten-new-1981}} (see also {\cite{parker-wittens-1982}}) which greatly
simplified the proof of the spacetime positive mass theorems when the initial
data set is spin. The spin condition is automatically satisfied in the
dimension 3 which is of more relevance to physics.

An initial data set $(M^n, g, p)$ is an $n$-dimensional manifold which arises
as a spacelike hypersurface of a Lorentzian manifold $(\mathcal{S}^{n, 1},
\tilde{g})$ with $p$ being the second fundamental form. The components $T_{0
0}$ and $T_{0 i}$ of the Einstein tensor (or the energy-momentum tensor) $T$
are respectively called the \text{{\itshape{energy}}}
\text{{\itshape{density}}} $\mu$ and the \text{{\itshape{current}}}
\text{{\itshape{density}}} $J$. Let $e_0$ be the future directed unit normal
of $M$ to $\mathcal{S}$, $e_i$ be an orthonormal basis of the tangent space of
$M$ and we use the convention on $p$ that $p_{i j} = \tilde{g}
(\tilde{\nabla}_{e_i} e_0, e_j)$.

The energy density by the Gauss equation is
\begin{equation}
  2 \mu = R_g + (\ensuremath{\operatorname{tr}}_g p)^2 - |p|_g^2
\end{equation}
and the current density by the Gauss-Codazzi equation is
\begin{equation}
  J =\ensuremath{\operatorname{div}}p - g \mathrm{d}
  (\ensuremath{\operatorname{tr}}_g p) .
\end{equation}
\begin{definition}
  We say that $(M, g, p)$ satisfies the interior dominant energy condition if
  \begin{equation}
    \mu \geqslant |J| . \label{interior dec}
  \end{equation}
  If $\partial M \neq \emptyset$, let $\eta$ be the outward normal of
  $\partial M$ in $M$, $H_{\partial M}
  =\ensuremath{\operatorname{div}}_{\partial M} \eta$. We say that $(M, g, p)$
  satisfies the tilted boundary dominant energy condition if
  \begin{equation}
    H_{\partial M} \pm \cos \theta \ensuremath{\operatorname{tr}}_{\partial M}
    p \geqslant \sin \theta |p (\eta, \cdot)^{\top} | \text{ on } \partial M,
    \label{tilt dec}
  \end{equation}
  where $\theta \in [0, \tfrac{\pi}{2}]$ is a constant angle and $p (\eta,
  \cdot)^{\top}$ denotes the component of the 1-form $p (\eta, \cdot)$
  tangential to $\partial M$.
\end{definition}

The tilted boundary dominant energy condition \eqref{tilt dec} generalizes the
tangential ($\theta = \pm \tfrac{\pi}{2}$) and normal boundary dominant energy
conditions ($\theta = 0$) in {\cite{almaraz-spacetime-2019}}. Now we recall
the definition of an asymptotically flat initial data set with a noncompact
boundary and its ADM energy and linear momentum of
{\cite{almaraz-spacetime-2019}}.

\begin{definition}
  We say that an initial data set $(M, g, p)$ is asymptotically flat with a
  noncompact boundary if there exists a compact set $K$ such that $M$ is
  diffeomorphic to the Euclidean half-space $\mathbb{R}^n_+ \backslash B_1$
  and
  \begin{equation}
    |g - \delta | + |x|  | \partial g| + |x|^2 | \partial^2 g| + |x| |p| +
    |x|^2 |p| = o (r^{- \tfrac{n - 2}{2}}), \label{decay}
  \end{equation}
  where $B_1$ is a standard Euclidean ball.
\end{definition}

\begin{definition}
  \label{ADM 4-vector}The quantities defined as

  \begin{equation}
    E = \lim_{r \to \infty} \left[ \int_{S^{n - 1, r}_+} (g_{i j, j} - g_{j j,
    i}) \nu^i - \int_{S^{n - 2, r}} e_{\alpha n} \vartheta^{\alpha} \right],
    \label{adm energy}
  \end{equation}
  and
  \begin{equation}
    P_i = 2 \int_{S^{n - 1, r}_+} \pi_{i j} \nu^j . \label{adm p}
  \end{equation}
  are respectively called the ADM energy and ADM linear momentum. Here, $\nu$
  is unit normal to $S_+^{n - 1, r}$, $\vartheta$ is normal to $S^{n - 2, r}$
  in $\partial M$ and $\pi = p  - g\ensuremath{\operatorname{tr}}_g p$. Denote
  $\hat{P} = (P_1, \ldots, P_{n - 1})$, $S^{n - 1, r}_+$ is the upper half of
  the coordinate sphere of radius $r$ and $S^{n - 2, r} = \partial S_+^{n - 1,
  r}$.
\end{definition}

Note that we include $P_n$ in the ADM linear-momentum. This is a difference
compared to {\cite{almaraz-spacetime-2019}}. We have the following two
spacetime positive mass theorems.

\begin{theorem}
  \label{tilt pmt}If $(M, g)$ is spin and $(M, g, p)$ satisfies the interior
  dominant energy condition \eqref{interior dec} and the tilted boundary
  dominant energy condition \eqref{tilt dec} for some nonzero $\theta \in (0,
  \tfrac{\pi}{2}]$, then
  \begin{equation}
    E \pm \cos \theta P_n \geqslant \sin \theta | \hat{P} | . \label{eq tilt
    pmt}
  \end{equation}
\end{theorem}

The special case $\theta = \tfrac{\pi}{2}$ of the theorem is due to
{\cite{almaraz-spacetime-2019}}. As we shall see later, \eqref{eq tilt pmt} is
related to an energy-momentum vector \eqref{tilt mass}, and the proof of
Theorem \ref{tilt pmt} already implies the particular case below.

\begin{theorem}
  \label{tangential pmt}If $M$ is spin, $\mu \geqslant |J|$ and $H \pm
  \ensuremath{\operatorname{tr}}_{\partial M} p \geqslant 0$, then
  \begin{equation}
    E \pm P_n \geqslant 0. \label{eq tangential pmt}
  \end{equation}
\end{theorem}

The time-symmetric case $p = 0$ of the theorem first appeared in
{\cite{almaraz-positive-2016}} where a minimal surface proof was also given.
The rough idea is: assume that the energy (mass) $E$ is negative, we can
perturb the metric so that it is harmonically flat, the scalar curvature and
the mean curvature of the boundary are strictly positive. Then the boundary
$\partial M$ and a plane asymptotically parallel to $\partial M$ serve as the
barriers and we can find an area-minimizing minimal surface which is
asymptotic to a coordinate plane that lies in between. Then the Gauss-Bonnet
theorem applied on the stable minimal plane contradicts the nonnegativity of
scalar curvature and mean curvature. An alternative proof was given by the
author {\cite{chai-positive-2018}}. Instead, \ the free boundary minimal
surface was used. Assume that $E < 0$, we are able to construct similarly a
free boundary area-minimizing surface that lies in between two coordinate
half-planes. The existence again contradicts the Gauss-Bonnet theorem.

Observing the two works, one should be able to conclude that the two proofs
using the minimal surface are actually proofs of two special cases when $p$
vanishes: (I) $\theta = \pi / 2$ in {\cite{chai-positive-2018}}; (II) or
$\theta = 0$ in {\cite{almaraz-positive-2016}}. This suggests that there is a
proof via stable minimal surface with capillary boundary conditions for the
case $p \equiv 0$ as well.

The capillarity also naturally arises in Gromov dihedral rigidity conjecture
{\cite{gromov-dirac-2014}}. The Euclidean version of Gromov dihedral rigidity
conjecture says that if a Riemannian polyhedron has nonnegative scalar
curvature, mean convex faces and its dihedral angles are less than its
Euclidean model, then it must be flat. The original motivation of Gromov
dihedral rigidity is a characterization of nonnegative scalar curvature in the
weak sense, it is also a localization of the positive mass theorems. Li
{\cite{li-polyhedron-2020}} confirmed this conjecture in some special cases
and the method he used is precisely minimal surface with a capillary boundary
condition. It is reasonable to establish directly Theorem \ref{tilt pmt}
implementing {\cite{eichmair-spacetime-2016}} using the capillary marginally
outer trapped surface {\cite{alaee-stable-2020}}. See Appendix \ref{capillary
MOTS} for a compact version of the positive mass theorem.

The article is organized as follows:

In Section \ref{sec:invariance}, we propose the mass related Theorems
\ref{tilt pmt} and \ref{tangential pmt} and show the invariance. In Section
\ref{sec dirac}, we collect basics of the chirality operator \eqref{Q} and the
hypersurface Dirac operator including the most important
Schrodinger-Lichnerowicz formula. In Section \ref{sec positive mass theorem},
we give the proofs of Theorems \ref{tilt pmt} and \ref{tangential pmt}. We
also include some partial results on the rigidity, that is, when equalities
are achieved in \eqref{eq tilt pmt} and \eqref{eq tangential pmt}. See
Propositions \ref{interior consequence} and \ref{boundary consequence}.

\section*{\text{{\bfseries{Acknowledgments}}}}

I would like to thank Xueyuan Wan (CQUT, Chongqing), Tin Yau Tsang (UCI) and
Martin Li (CUHK) for various discussions, and Levi Lima (UFC, Brazil) for
communicating the preprint {\cite{almaraz-rigidity-2022-arxiv}} which puts
Theorem \ref{tangential pmt} in a more general context of Theorem \ref{tilt
pmt}. Part of this work was carried out when I was a member of Korea Institute
for Advanced Study (KIAS) under the grant No. MG074402. I am also partially
supported by National Research Foundation of Korea grant No. 2002R1C1C1013511.

\section{The invariance of mass}\label{sec:invariance}

In this section, we introduce the energy-momentum vector $(E^{\theta},
P_i^{\theta})$ in \eqref{tilt mass} based on the Hamiltonian analysis
performed in {\cite{hawking-gravitational-1996}} and point out that the tilted
dominant energy condition \eqref{tilt dec} appears in selecting a suitable
lapse function and the shift vector.

\subsection{Hamiltonian formulation and mass invariance}

Assume at present that $M$ is compact, we infinitesimally deform the initial
data set $(M, g, p)$ in $\mathcal{S}^{n, 1}$ in the direction of a future
directed timelike vector field $T$. Let $\phi_s$ be the local flow of $T$,
$M_s = \phi_s (M)$. We assume that the unit normal $e_0$ to $M_s$ is always
tangential to the timelike hypersurface foliated by $\partial M_s$. Let $T = N
e_0 + X$, where $N$ is called the \text{{\itshape{lapse}}} function and the
vector field $X$ tangent to $M$ is called the \text{{\itshape{shift vector}}},
then the Hamiltonian along $M$ is given by (see
{\cite{hawking-gravitational-1996}})
\begin{equation}
  \mathcal{H} (N, X) = \int_M [N \mu + 2 J (X)] + 2 \int_{\partial M} [N H - p
  (X, \eta) +\ensuremath{\operatorname{tr}}_g p \langle X, \eta \rangle] .
  \label{hamiltonian}
\end{equation}
The tilted boundary dominant energy condition \eqref{tilt dec} now comes from
selecting $N = 1$ and $X = \cos \theta \eta + \sin \theta \tau$ where $\tau$
is tangent to $\partial M$ in the boundary term of the Hamiltonian. Indeed,
\begin{align}
& N H - p (X, \eta) +\ensuremath{\operatorname{tr}}_g p \langle X, \eta
\rangle \\
= & H - \cos \theta p (\eta, \eta) - \sin \theta p (\tau, \eta) + \cos
\theta \ensuremath{\operatorname{tr}}_g p \\
= & H + \cos \theta \ensuremath{\operatorname{tr}}_{\partial M} p - \sin
\theta p (\tau, \eta) \\
\geqslant & 0
\end{align}
if \eqref{tilt dec} holds.

Now let $(M, g, p)$ be the background $(\mathbb{R}^n_+, \delta, 0)$, we take
$N$ to be a constant and $X$ be a translational Killing vector field of
$\mathbb{R}^n_+$. We consider the Hamiltonian $\mathcal{H}_{\varepsilon} (N,
X)$ on $(M, \delta + \varepsilon g, \varepsilon p)$ with $(g, p)$ satisfying
\eqref{decay}. We do the Taylor expansion of $\mathcal{H}_{\varepsilon}$ with
respect to $\varepsilon$, due to the fact that $M$ is noncompact, usually the
first order terms do not vanish. These terms evaluated at infinity are
precisely those given in Definition \ref{ADM 4-vector}. For a more complete
account of these facts, we refer the readers to
{\cite{hawking-gravitational-1996}}, {\cite{michel-geometric-2011}} and
{\cite{almaraz-spacetime-2019}}.

We define the \text{{\itshape{charge}}} \text{{\itshape{density}}} which is a
1-form,
\begin{align}
& \mathbb{U}_{(g, k)} (N, X) \\
= & N (\ensuremath{\operatorname{div}}_{\delta} g - \mathrm{d}
(\ensuremath{\operatorname{tr}}_{\delta} g)) - (g - \delta) (\nabla^{\delta}
N, \cdot) \\
& \quad +\ensuremath{\operatorname{tr}}_{\delta} (g - \delta) \mathrm{d} N
+ 2 (p (X, \cdot) -\ensuremath{\operatorname{tr}}_{\delta} p \langle \cdot,
X \rangle_{\delta}) .
\end{align}
Let $\mathcal{T}$ be the space of translational Killing vector fields of
Minkowski spacetime denoted by $\mathbb{R}^{1, n}$. It is easy to see that
$\mathcal{T}$ is identified with $\mathbb{R} \oplus W$ with $\mathbb{R}$
factor representing the translation in a chosen timelike direction
$\partial_0$ and $W$ being the linear space spanned by all translational
Killing vector fields of $(\mathbb{R}^n, \delta)$ orthogonal to $\partial_0$.
Each $T \in \mathcal{T}$ can be uniquely written in the form $T = N \partial_0
+ X^i \partial_i$ where $N \in \mathbb{R}$ and $X^i \in \mathbb{R}$. We define
the energy-momentum functional as follows:
\begin{equation}
  \mathcal{M} (T) = \lim_{r \to \infty} \left[ \int_{S^{n - 1, r}_+}
  \mathbb{U}_{(g, k)} (N, X) + \int_{S^{n - 2, r}} N g (\bar{\eta},
  \bar{\vartheta}) \right] .
\end{equation}
It was shown in {\cite[Proposition 3.3]{almaraz-spacetime-2019}} that the
energy-momentum functional $\mathcal{M} (T)$ does not depend on the asymptotic
coordinates (fixing $\partial_0$) chosen at infinity.

For any $\theta \in (0, \tfrac{\pi}{2}]$, we define
\begin{equation}
  E^{\theta} =\mathcal{M} (\tfrac{1}{\sin \theta} \partial_0 + \tfrac{\cos
  \theta}{\sin \theta} \partial_n), P_i^{\theta} =\mathcal{M} (\partial_i)
  \text{ for any } i \neq n. \label{tilt mass}
\end{equation}
It is easy to check that $E =\mathcal{M} (\partial_0)$, $P_i =\mathcal{M}
(\partial_i)$ where $(E, P)$ is as defined in Definition \ref{ADM 4-vector},
so $E^{\theta} = \tfrac{1}{\sin \theta} E + \tfrac{\cos \theta}{\sin \theta}
P_n$. We have the following.

\begin{theorem}
  Given any asymptotically flat initial data set $(M, g, k)$, for any $\theta
  \in (0, \tfrac{\pi}{2}]$, the vector $(E^{\theta}, P^{\theta}) \in
  \mathbb{R}^{1, n - 1}$ is well defined (up to composition with an element of
  $S O_{1, n - 1}$). In particular,
  \begin{equation}
    - (E^{\theta})^2 + \sum_{i \neq n} (P_i^{\theta})^2
  \end{equation}
  and the causal character $(E^{\theta}, P^{\theta}) \in \mathbb{R}^{1, n -
  1}$ do not depend on the chart at infinity to compute $(E^{\theta},
  P^{\theta})$.
\end{theorem}

\begin{proof}
  Let $\tilde{\partial}_0 = \tfrac{1}{\sin \theta} (\partial_0 + \cos \theta
  \partial_n)$, $\tilde{\partial}_n = \tfrac{1}{\sin \theta} (\cos \theta
  \partial_0 + \partial_n)$ and $\tilde{\partial}_i = \partial_i$. There is a
  Lorentz boost from $(\partial_0, \partial_1, \ldots, \partial_n)$ to
  $(\tilde{\partial}_0, \tilde{\partial}_1, \ldots, \tilde{\partial}_n)$ such
  that
  \begin{equation}
    \left(\begin{array}{c}
      \tilde{\partial}_0\\
      \tilde{\partial}_n
    \end{array}\right) = \left(\begin{array}{cc}
      \cosh \rho & \sinh \rho\\
      \sinh \rho & \cosh \rho
    \end{array}\right) \left(\begin{array}{c}
      \partial_0\\
      \partial_n
    \end{array}\right),
  \end{equation}
  on the plane spanned by $\{\partial_0, \partial_n \}$ with $\rho$ defined by
  $\cosh \rho = \tfrac{1}{\sin \theta}$. So $(\tilde{\partial}_0,
  \tilde{\partial}_1, \ldots, \tilde{\partial}_n)$ gives a new coordinate
  system for the Minkowski spacetime $\mathbb{R}^{1, n}$. Let $(\tilde{x}_0,
  \tilde{x}_1, \ldots, \tilde{x}_n) \in \mathbb{R}^{1, n}$ where $\tilde{x}$
  is expressed in the new coordinates. Obviously,
  \begin{equation}
    - (\tilde{x}_0)^2 + \sum_{i \neq n} (\tilde{x}_n)^2
  \end{equation}
  is invariant under linear Lorentz transformations of $\mathbb{R}^{1, n}$
  which fixes $\tilde{\partial}_n$. These transformations as a subgroup of the
  special Lorentz group $S O_{1, n}$ is isomorphic to $S O_{1, n - 1}$. The
  discussion applies to
  \[ (\mathcal{M}(\tilde{\partial}_0), \mathcal{M}(\widetilde{\partial_1}),
     \ldots, \mathcal{M}(\tilde{\partial}_{n - 1}),
     \mathcal{M}(\tilde{\partial}_n)), \]
  and this is our theorem.
\end{proof}

For the case $\theta = 0$, it is simpler.

\begin{theorem}
  Given any asymptotically flat initial data set $(M, g, k)$, the quantity $E
  \pm P_n$ is a numerical invariant under isometries of $\mathbb{R}^n_+$ which
  includes rotations and translations of the $(n - 1)$-dimensional hyperplane
  $\partial \mathbb{R}^n_+$.
\end{theorem}

\begin{proof}
  Note that $E$ and $P_n$ are invariant under rotations and translations of
  the hyperplane $\{\partial_1, \ldots, \partial_{n - 1} \}$, see
  {\cite[Proposition 3.3]{almaraz-spacetime-2019}}.
\end{proof}

\section{Hypersurface Dirac operator}\label{sec dirac}

In this section, we recall the hypersurface Dirac spinors and the related
Schrodinger-Lichnerowicz formula \eqref{integral slw}. We review the chirality
operator \eqref{Q} and we relate the boundary condition \eqref{boundary
condition} to the geometric quantities along the boundary $\partial M$ in
Lemma \ref{boundary calculation lemma}.

\subsection{Hypersurface Dirac operator}

The standard reference of spin geometry is {\cite{lawson-spin-1989}}, we also
refer to {\cite{parker-wittens-1982}}, {\cite{hijazi-dirac-witten-2003}}.
Denote by $\mathbb{S}$ the local spinor bundle of $\mathcal{S}$, since $M$ is
spin, $\mathbb{S}$ exists globally over $M$. This spinor bundle $\mathbb{S}$
is called the hypersurface spinor bundle of $M$. Let $\tilde{\nabla}$ and
$\nabla$ denote respectively the Levi-Civita connections of $\tilde{g}$ and
$g$, we use the same symbols to denote the lifts of the connections to the
hypersurface spinor bundle.

There exists a Hermitian inner product $(\cdot, \cdot)$ on $\mathbb{S}$ over
$M$ which is compatible with the spin connection $\tilde{\nabla}$. For any
1-form $\omega$ of $\mathcal{S}$ and the hypersurface spinors $\phi$, $\psi$,
we have
\[ (\omega \cdot \phi, \psi) = (\phi, \omega \cdot \psi) \]
where the dot $\cdot$ denotes the Clifford multiplication. This inner product
is not positive definite. However, there exists on $\mathbb{S}$ over $M$ a
positive definite Hermitian inner product defined by
\[ \langle \phi, \psi \rangle = (e^0 \cdot \phi, \psi) \]
where $e^0$ is the future-directed unit timelike normal to $M$. We see that
\begin{equation}
  \langle e^0 \cdot \phi, \psi \rangle = \langle \phi, e^0 \cdot \psi \rangle,
  \text{ } \langle e^i \cdot \phi, \psi \rangle = - \langle \phi, e^i \cdot
  \psi \rangle,
\end{equation}
where $\{e^i \}$ are the dual frame of an orthonormal basis $\{e_i \}$ over
$M$. Then the spinor connection $\tilde{\nabla}$ over $\mathbb{S}$ is related
to $\nabla$ by
\begin{equation}
  \tilde{\nabla}_i = \nabla_i - \tfrac{1}{2} p_{i j} e^j \cdot e^0 \cdot .
\end{equation}
This is essentially the spinorial Gauss equation. Moreover, the connection
$\nabla$ is compatible with $\langle \cdot, \cdot \rangle$ and $\nabla_i (e^0
\cdot \phi) = e^0 \cdot \nabla_i \phi$.

The hypersurface Dirac (or Dirac-Witten) operator is then defined by
\begin{equation}
  \tilde{D} = e^i \cdot \tilde{\nabla}_i = D + \tfrac{1}{2}
  \ensuremath{\operatorname{tr}}_g p e^0 \cdot,
\end{equation}
where $D$ is the standard Dirac operator. We call a spinor $\phi$ satisfying
$\tilde{D} \phi = 0$ a \text{{\itshape{(spacetime)}}}
\text{{\itshape{harmonic}}} \text{{\itshape{spinor}}}.

When our spacetime $\mathcal{S}$ is of dimension $3 + 1$, the local spinor
bundle $\mathbb{S}$ have a simpler algebraic description by the representation
theory of the special linear group $\ensuremath{\operatorname{SL}} (2,
\mathbb{C})$. In this case, the theory is easier to understand, see
{\cite[Section 2]{parker-wittens-1982}}.

We now collect relevant facts regarding $\tilde{\nabla}$ and $\tilde{D}$.

\begin{lemma}[{\cite{parker-wittens-1982}}]
  \label{adjoints}The adjoint of $\tilde{\nabla}$ is given by
\begin{align}
\tilde{\nabla}_i^{\ast} \psi & = - \nabla_i \psi - \tfrac{1}{2} p_{i j}
e^j \cdot e^0 \psi, \\
\mathrm{d} (\langle \phi, \psi \rangle \ast e^i) & = [\langle
\tilde{\nabla}_i \phi, \psi \rangle - \langle \phi,
\tilde{\nabla}^{\ast}_i \psi \rangle] \ast 1.
\end{align}
  Here $\ast$ denotes the Hodge operator. The Dirac operator $\tilde{D}$ is
  self-adjoint with
  \begin{equation}
    \mathrm{d} (\langle e^i \cdot \phi, \psi \rangle \ast e^i) = (\langle
    \tilde{D} \phi, \psi \rangle - \langle \phi, \tilde{D} \psi \rangle) \ast
    1.
  \end{equation}
  The Schrodinger-Lichnerowicz formula is given by
  \begin{equation}
    \tilde{D}^2 - \tilde{\nabla}^{\ast} \tilde{\nabla} = \tfrac{1}{2} (\mu - J
    \cdot e^0 \cdot) \label{slw} .
  \end{equation}
\end{lemma}

The integration form of the Schrodinger-Lichnerowicz formula is a direct
corollary of Lemma \ref{adjoints}.

\begin{theorem}[{\cite{parker-wittens-1982}}]
  Let $\Omega$ be a compact manifold with boundary, we have for any smooth
  spinor $\phi$ that
\begin{align}
& \int_{\Omega} | \tilde{D} \phi |^2 - \int_{\Omega} | \tilde{\nabla}
\phi |^2 + \int_{\partial \Omega} [\langle e^i \cdot \tilde{D} \phi, \phi
\rangle + \langle \phi, \tilde{\nabla}_i \phi \rangle] \ast e^i
\\
= & \tfrac{1}{2} \int_{\Omega} \langle (\mu - J \cdot e^0 \cdot) \phi,
\phi \rangle . \label{integral slw}
\end{align}
\end{theorem}

\subsection{Boundary chirality operator}

We fix the convention. We use the Greek letters $\alpha$, $\beta$, $\gamma$ to
indicate the indices which are not $n$ in the rest of the paper. The vector
$e_n$ is used to denote the outer normal of $\partial M$ in $M$ and $h$
denotes the the second fundamental form of $\partial M$ given by $h_{\alpha
\beta} = \langle \nabla_{e_{\alpha}} e_n, e_{\beta} \rangle$, then $H : =
H_{\partial M} = \sum_{\alpha} h_{\alpha \alpha}$.

The following chirality operator $Q$ was introduced by {\cite[Definition
3.3]{almaraz-rigidity-2022-arxiv}} where $e^0 \cdot e^n \cdot$ is replaced by
the Clifford multiplication of the complex volume element.

\begin{definition}
  Let $\theta \in [- \tfrac{\pi}{2}, \tfrac{\pi}{2}]$, define $Q$ by
  \begin{equation}
    Q  \phi = \cos \theta e^0 \cdot e^n \cdot \phi + \sqrt{- 1} \sin \theta
    e^n \cdot \phi . \label{Q}
  \end{equation}
\end{definition}

In \eqref{Q}, we require that $\theta \in [- \tfrac{\pi}{2}, \tfrac{\pi}{2}]$,
but in \eqref{tilt dec} we require that $\theta \in [0, \tfrac{\pi}{2}]$. The
reason is that we need a choice on the sign of $\theta$ later in the proof of
Theorem \ref{tilt pmt}. We collect a few commutative and anti-commutative
properties of $Q$ below.

\begin{lemma}
  \label{Q algebraic}For the operator $Q$, we have
  \begin{enumerate}
    \item $Q^2 = 1$ and $Q$ is self-adjoint;
    
    \item $e^n \cdot Q + Q \cdot e^n = - 2 \sqrt{- 1} \sin \theta$;
    
    \item $e^{\alpha} \cdot e^{\beta} \cdot e^n \cdot Q + Q e^{\alpha} \cdot
    e^{\beta} \cdot e^n \cdot = - 2 \sqrt{- 1} \sin \theta e^{\alpha} \cdot
    e^{\beta} \cdot$;
    
    \item $e^0 \cdot Q + Q \cdot e^0 = 0$;
    
    \item $e^{\alpha} \cdot e^{\beta} \cdot e^0 \cdot Q + Q e^{\alpha} \cdot
    e^{\beta} \cdot e^0 \cdot = 0$;
    
    \item $e^{\alpha} \cdot Q \phi - Q (e_{\alpha} \cdot \phi) = 2 \sqrt{- 1}
    \sin \theta e^{\alpha} \cdot e^n \cdot \phi$;
    
    \item $e^{\alpha} \cdot Q \phi + Q (e_{\alpha} \cdot \phi) = 2 \cos \theta
    e^{\alpha} \cdot e^0 \cdot e^n \cdot \phi$
    
    \item $e^{\alpha} \cdot e^n \cdot Q + Q \cdot e^{\alpha} \cdot e^n = 0$;
    
    \item $e^n \cdot e^0 \cdot Q + Q e^n \cdot e^0 = - 2 \cos \theta$;
    
    \item $e^{\alpha} \cdot e^{\beta} \cdot e^0 \cdot e^n \cdot Q + Q
    e^{\alpha} \cdot e^{\beta} \cdot e^0 \cdot e^n \cdot = 2 \cos \theta
    e^{\alpha} \cdot e^{\beta} \cdot$.
    
    \item $e^{\alpha} \cdot e^0 \cdot Q + Q e^{\alpha} \cdot e^0 = 2 \sqrt{-
    1} \sin \theta e^{\alpha} \cdot e^0 \cdot e^n$;
    
    \item For $\alpha \neq \beta$, $e^{\alpha} \cdot e^{\beta} \cdot Q = Q
    e^{\alpha} \cdot e^{\beta} \cdot$;
    
    \item $e^{\alpha} \cdot e^n \cdot e^0 \cdot Q = Q \cdot e^{\alpha} \cdot
    e^n \cdot e^0$.
  \end{enumerate}
\end{lemma}

\begin{proof}
  All the items follows from direct calculation starting from the definition
  of $Q$ in \eqref{Q}. As an example, we only show the last item. By \eqref{Q}
  and anti-commutative property of the Clifford multiplication,
\begin{align}
e^{\alpha} \cdot e^n \cdot e^0 \cdot Q \phi & = \cos \theta e^{\alpha}
\cdot e^n \cdot e^0 \cdot e^0 \cdot e^n \cdot \phi + \sqrt{- 1} \sin
\theta e^{\alpha} \cdot e^n \cdot e^0 \cdot e^n \cdot \phi \\
& = - \cos \theta e^{\alpha} \cdot \phi + \sqrt{- 1} \sin \theta
e^{\alpha} \cdot e^0 \cdot \phi, \\
Q (e_{\alpha} \cdot e^n \cdot e^0 \cdot \phi) & = \cos \theta e^0 \cdot
e^n \cdot e^{\alpha} \cdot e^n \cdot e^0 \cdot \phi + \sqrt{- 1} \sin
\theta e^n \cdot e^{\alpha} \cdot e^n \cdot e^0 \cdot \phi \\
& = - \cos \theta e^{\alpha} \cdot \phi + \sqrt{- 1} \sin \theta
e^{\alpha} \cdot e^0 \cdot \phi .
\end{align}
  So we know that the last item holds true.
\end{proof}

\subsection{Boundary terms in Schrodinger-Lichnerowicz formula}

We calculate
\begin{equation}
  [\langle e^i \cdot \tilde{D} \phi, \phi \rangle + \langle \phi,
  \tilde{\nabla}_i \phi \rangle] \ast e^i
\end{equation}
along $\partial M$ when
\begin{equation}
  Q \phi = \pm \phi . \label{boundary condition}
\end{equation}
First, we compute a few inner products of spinors satisfying \eqref{boundary
condition}.

\begin{lemma}
  If a spinor $\phi$ satisfies \eqref{boundary condition} along $\partial M$,
  then
\begin{align}
\langle \sqrt{- 1} e^n \cdot \phi, \phi \rangle & = \pm \sin \theta | \phi
|^2, \label{en directional norm} \\
\langle e^n \cdot e^0 \cdot \phi, \phi \rangle & = \mp \cos \theta | \phi
|^2, \label{en e0 component} \\
\langle e^{\alpha} \cdot e^0 \cdot \phi, \phi \rangle & = \mp \langle
\sqrt{- 1} \sin \theta e^{\alpha} \cdot e^n \cdot e^0 \cdot \phi, \phi
\rangle . \label{ej e0 component}
\end{align}
\end{lemma}

\begin{proof}
  The first term \eqref{en directional norm} already appeared in
  {\cite[Proposition 3.11]{almaraz-rigidity-2022-arxiv}}. From Lemma \ref{Q
  algebraic}, we have
  \[ \langle \sqrt{- 1} e^n \cdot Q \phi, \phi \rangle + \langle Q \cdot
     \sqrt{- 1} e^n \cdot \phi, \phi \rangle = 2 \sin \theta | \phi |^2 . \]
  Since $Q$ is self-adjoint, so
  \begin{equation}
    \langle \sqrt{- 1} e^n \cdot Q \phi, \phi \rangle + \langle \sqrt{- 1} e^n
    \cdot \phi, Q \phi \rangle = 2 \sin \theta | \phi |^2 .
  \end{equation}
  Because $Q \phi = \pm \phi$, we have
  \begin{equation}
    \pm 2 \langle \sqrt{- 1} e^n \cdot \phi, \phi \rangle = 2 \sin \theta |
    \phi |^2,
  \end{equation}
  which is the first item. The rest follow similarly from corresponding
  relations from Lemma \ref{Q algebraic}.
\end{proof}

The following lemma relates the boundary term in the integration form of
Schrodinger-Lichnerowicz formula \eqref{integral slw} with the mean curvature
$H$, $\ensuremath{\operatorname{tr}}_{\partial M} p$, $p_{n j}$ along the
boundary, and in particular, the tilted boundary dominant energy condition
\eqref{tilt dec}.

\begin{lemma}
  \label{boundary calculation lemma}If a spinor $\phi$ satisfies
  \eqref{boundary condition} along $\partial M$, then
\begin{align}
& \langle \tilde{\nabla}_{e^n} \phi + e^n \cdot \tilde{D} \phi, \phi
\rangle \\
= & \langle D^{\partial M} \phi, \phi \rangle - \tfrac{1}{2} H | \phi |^2
\mp \tfrac{1}{2} \cos \theta \ensuremath{\operatorname{tr}}_{\partial M} p
| \phi |^2 \\
& \text{ } \pm \tfrac{1}{2} \sin \theta \langle \sqrt{- 1} p_{n \gamma}
e^{\gamma} \cdot e^n \cdot e^0 \cdot \phi, \phi \rangle \label{boundary
calculation} .
\end{align}
\end{lemma}

\begin{proof}
  Let $D^{\partial M}$ be a boundary Dirac operator{\underline{}}or defined by
  \[ D^{\partial M} = e^n \cdot e^{\alpha} \cdot \nabla_{\alpha}^{\partial M}
     . \]
  Here, $\nabla^{\partial M}$ is the spin connection intrinsic to $\partial M$
  explicitly defined on spinor fields on $M$ restricted to $\partial M$ as
  \[ \nabla^{\partial M}_{\alpha} = \nabla_{\alpha} - \tfrac{1}{2} h_{\alpha
     \beta} e^n \cdot e^{\beta} \cdot . \]
  We calculate $D^{\partial M} \phi$ with $\phi$ satisfying \eqref{boundary
  condition} and we see
\begin{align}
& D^{\partial M} \phi \\
= & e^n \cdot e^{\alpha} \cdot (\nabla_{\alpha} \phi - \tfrac{1}{2}
h_{\alpha \beta} e^n \cdot e^{\beta} \cdot \phi) \\
= & e^n \cdot (D \phi - e^n \cdot \nabla_{e^n} \phi) + \tfrac{1}{2} H \phi
\\
= & e^n \cdot D \phi + \nabla_{e^n} \phi + \tfrac{1}{2} H \phi \\
= & e^n \cdot (\tilde{D} \phi - \tfrac{1}{2}
\ensuremath{\operatorname{tr}}_g p e^0 \cdot \phi) + (\tilde{\nabla}_{e^n}
\phi + \tfrac{1}{2} p_{n j} e^j \cdot e^0 \cdot \phi) + \tfrac{1}{2} H
\phi
\end{align}
  So we have that
\begin{align}
& \langle \tilde{\nabla}_{e^n} \phi + e^n \cdot \tilde{D} \phi, \phi
\rangle \\
= & \langle D^{\partial M} \phi, \phi \rangle + \langle \tfrac{1}{2}
\ensuremath{\operatorname{tr}}_g p e^n \cdot e^0 \cdot \phi - \tfrac{1}{2}
p_{n j} e^j \cdot e^0 \cdot \phi, \phi \rangle - \tfrac{1}{2} H | \phi |^2
.
\end{align}
  From \eqref{en e0 component}, we have
\begin{align}
& \left\langle \tfrac{1}{2} \ensuremath{\operatorname{tr}}_g p e^n \cdot
e^0 \cdot \phi - \tfrac{1}{2} p_{n j} e^j \cdot e^0 \cdot \phi, \phi
\right\rangle \\
= & \left\langle \tfrac{1}{2} (\ensuremath{\operatorname{tr}}_{\partial M}
p e^n + p_{n n} e^n) \cdot e^0 \cdot \phi - \tfrac{1}{2} (p_{n \alpha}
e^{\alpha} + p_{n n} e^n) \cdot e^0 \cdot \phi, \phi \right\rangle
\\
= & \tfrac{1}{2} \ensuremath{\operatorname{tr}}_{\partial M} p \langle e^n
\cdot e^0 \cdot \phi, \phi \rangle - \tfrac{1}{2} \langle p_{n \alpha}
e^{\alpha} \cdot e^0 \cdot \phi, \phi \rangle \\
= & \mp \tfrac{1}{2} \cos \theta \ensuremath{\operatorname{tr}}_{\partial
M} p | \phi |^2 - \tfrac{1}{2} \langle p_{n \alpha} e^{\alpha} \cdot e^0
\cdot \phi, \phi \rangle \\
= & \mp \tfrac{1}{2} \cos \theta \ensuremath{\operatorname{tr}}_{\partial
M} p | \phi |^2 \mp \tfrac{1}{2} \sin \theta \langle \sqrt{- 1} p_{n
\alpha} e^{\alpha} \cdot e^0 \cdot e^n \cdot \phi, \phi \rangle
\end{align}
  which follows from \eqref{en e0 component} and \eqref{ej e0 component}.
\end{proof}

\section{the positive mass theorem}\label{sec positive mass theorem}

With the help of results from the previous sections, in this section, we
finish the proofs of Theorems \ref{tilt pmt} and \ref{tangential pmt}. We give
some consequences of vanishing mass, that is,
\begin{equation}
  E \pm \cos \theta P_n = \sin \theta | \hat{P} | . \label{vanishing mass}
\end{equation}
See Propositions \ref{interior consequence} and \ref{boundary consequence}.

\subsection{Existence of a spacetime harmonic spinor}

When the initial data set $(M, g, p)$ is flat and totally geodesic, i.e. $(M,
g, p)$ is $(\mathbb{R}^n_+, \delta, 0)$, we define
\begin{equation}
  \bar{Q} \phi = \cos \theta \mathrm{d} x^0 \cdot \mathrm{d} x^n \cdot \phi +
  \sqrt{- 1} \sin \theta \mathrm{d} x^n \cdot \phi .
\end{equation}
Note that $\bar{Q}^2 = I$, the eigenvalues of $\bar{Q}$ are $\pm 1$. The
standard hypersurface spinor bundle $\bar{\mathbb{S}}$ over $(\mathbb{R}_+^n,
\delta, 0)$ splits into two eigen subbundles and the spinor $\phi$ satisfying
\begin{equation}
  \bar{Q} \phi = \pm \phi \label{background eigen spinor}
\end{equation}
is closely related to our problem.

We recall the following existence of a spacetime harmonic spinor $\phi$ which
is asymptotic to a constant spinor $\phi_0$ satisfying \eqref{background eigen
spinor}, and we extract the mass from the boundary integral in \eqref{integral
slw}. By {\cite[Proposition 5.3]{almaraz-spacetime-2019}} and the discussions
that followed, we have the following.

\begin{theorem}
  \label{existence of hypersurface dirac operator}Assume that $(M, g, k)$
  satisfies the dominant energy conditions \eqref{interior dec} and
  \eqref{tilt dec}, given any nonzero constant spinor $\phi_0$ satisfying
  \eqref{background eigen spinor} , there exists a spinor $\phi$ which is
  asymptotic to $\phi_0$ and satisfies
\begin{align}
\tilde{D} \phi & = 0 \text{ in } M, \\
Q \phi & = \pm \phi \text{ on } \partial M.
\end{align}
\end{theorem}

\subsection{Proof of positive mass theorems}

Using the $\phi$ of Theorem \ref{existence of hypersurface dirac operator} in
\eqref{integral slw}, we can give a proof of Theorem \ref{tilt dec} and the
proof works equally well for Theorem \ref{tangential pmt}.

\begin{proof}[Proof of Theorem \ref{tilt pmt}]
  Let $M_r$ be the compact region bounded by $\partial M$ and $S_r^{n - 1,
  +}$. By the integral form of Schrodinger-Lichnerowicz formula
  \eqref{integral slw}, we have for any spinor $\phi$, we have
\begin{align}
& \int_{M_r} | \tilde{D} \phi |^2 - \int_{M_r} | \tilde{\nabla} \phi |^2
+ \int_{\partial M_r} [\langle e^i \cdot \tilde{D} \phi, \phi \rangle +
\langle \phi, \tilde{\nabla}_i \phi \rangle] \ast e^i \\
= & \tfrac{1}{2} \int_{M_r} \langle (\mu - J \cdot e^0 \cdot) \phi, \phi
\rangle .
\end{align}
  Note that $\partial M_r$ are made of two portions: one lies in the interior
  of $M$ and the other lies on $\partial M$. We require that $Q \phi = \pm
  \phi$ along $\partial M$, so by Lemma \ref{boundary calculation lemma},
\begin{align}
& \int_{M_r} | \tilde{D} \phi |^2 - \int_{M_r} | \tilde{\nabla} \phi |^2
+ \int_{\partial M_r \cap \ensuremath{\operatorname{int}}M} [\langle e^i
\cdot \tilde{D} \phi, \phi \rangle + \langle \phi, \tilde{\nabla}_i \phi
\rangle] \ast e^i \\
& + \int_{\partial M_r \cap \partial M} \langle D^{\partial M} \phi, \phi
\rangle - \tfrac{1}{2} H | \phi |^2 \mp \tfrac{1}{2} \cos \theta
\ensuremath{\operatorname{tr}}_{\partial M} p | \phi |^2 \\
& \quad \pm \int_{\partial M_r \cap \partial M} \tfrac{1}{2} \sin \theta
\langle \sqrt{- 1} p_{n \gamma} e^{\gamma} \cdot e^n \cdot e^0 \cdot \phi,
\phi \rangle \\
= & \tfrac{1}{2} \int_{M_r} \langle (\mu - J \cdot e^0 \cdot) \phi, \phi
\rangle .
\end{align}
  It follows that $\langle D^{\partial M} \phi, \phi \rangle = 0$ from
  {\cite[(4.27)]{chrusciel-mass-2003}} (with $\varepsilon$ there replaced by
  $Q$). We have
\begin{align}
& \int_{\partial M_r \cap \ensuremath{\operatorname{int}}M} [\langle e^i
\cdot \tilde{D} \phi, \phi \rangle + \langle \phi, \tilde{\nabla}_i \phi
\rangle] \ast e^i \\
\to & \tfrac{1}{4} (E \pm \cos \theta P_n) | \phi |^2 \pm \sin \theta
P_{\gamma} \langle \sqrt{- 1} \mathrm{d} x^{\gamma} \cdot \mathrm{d} x^n
\cdot \mathrm{d} x^0 \cdot \phi_0, \phi_0 \rangle_{\delta} \label{mass}
\end{align}
  as $r \to \infty$. We calculate
\begin{align}
& [\langle e^i \cdot \tilde{D} \phi, \phi \rangle + \langle \phi,
\tilde{\nabla}_i \phi \rangle] \ast e^i \\
= & \langle \nu^{\flat} \cdot \tilde{D} \phi, \phi \rangle + \langle \phi,
\tilde{\nabla}_{\nu} \phi \rangle \\
= & \langle \nu^{\flat} \cdot D \phi, \phi \rangle + \langle \phi,
\nabla_{\nu} \phi \rangle \\
& + \tfrac{1}{2} (\ensuremath{\operatorname{tr}}_g p \nu^i - p_{i j}
\nu^j) \langle e^i \cdot e^0 \cdot \phi, \phi \rangle \\
= & \langle \nu^{\flat} \cdot D \phi, \phi \rangle + \langle \phi,
\nabla_{\nu} \phi \rangle - \tfrac{1}{2} \pi_{i j} \nu^j \langle e^i \cdot
e^0 \cdot \phi, \phi \rangle,
\end{align}
  where $\nu$ is the unit normal of $\partial M_r \cap
  \ensuremath{\operatorname{int}}M$ pointing to the infinity. Because that
  $\phi$ converges to a constant spinor $\phi_0$, we know that as $r \to
  \infty$
  \begin{equation}
    \int_{\partial M_r \cap \ensuremath{\operatorname{int}}M} \langle
    \nu^{\flat} \cdot D \phi, \phi \rangle + \langle \phi, \nabla_{\nu} \phi
    \rangle \to \tfrac{1}{4} E | \phi_0 |^2_{\delta}
  \end{equation}
  from {\cite[Section 5.2]{almaraz-positive-2016}} and

  \begin{equation}
    - \tfrac{1}{2} \int_{\partial M_r \cap \ensuremath{\operatorname{int}}M}
    \pi_{i j} \nu^j \langle e^i \cdot e^0 \cdot \phi, \phi \rangle \to -
    \tfrac{1}{4} P_i \langle \mathrm{d} x^i \cdot \mathrm{d} x^0 \cdot \phi_0,
    \phi_0 \rangle_{\delta}
  \end{equation}
  since $\langle e^i \cdot e^0 \cdot \phi, \phi \rangle$ converges to the
  constant $\langle \mathrm{d} x^i \cdot \mathrm{d} x^0 \cdot \phi_0, \phi_0
  \rangle_{\delta}$. Here $\delta$ is the standard Euclidean metric. From
  \eqref{en e0 component} that
  \begin{equation}
    P_n \langle \mathrm{d} x^n \cdot \mathrm{d} x^0 \cdot \phi_0, \phi_0
    \rangle_{\delta} = \mp \cos \theta P_n | \phi_0 |^2_{\delta},
  \end{equation}
  and from \eqref{ej e0 component} that
  \begin{equation}
    P_{\gamma} \langle \mathrm{d} x^{\gamma} \cdot \mathrm{d} x^0 \cdot
    \phi_0, \phi_0 \rangle_{\delta} = \mp \sin \theta P_{\gamma} \langle
    \sqrt{- 1} \mathrm{d} x^{\gamma} \cdot \mathrm{d} x^n \cdot \mathrm{d} x^0
    \cdot \phi_0, \phi_0 \rangle_{\delta} . \label{other pi}
  \end{equation}
  By the lemma below, we can make a choice of the constant spinor $\phi_0$ and
  the sign of $\theta$ such that
  \begin{equation}
    P_{\gamma} \langle \mathrm{d} x^{\gamma} \cdot \mathrm{d} x^0 \cdot
    \phi_0, \phi_0 \rangle_{\delta} = \sin | \theta |  | \hat{P} |  | \phi_0
    |_{\delta} .
  \end{equation}
  So taking limits as $r \to \infty$, we have
\begin{align}
& \tfrac{1}{4} (E \pm \cos \theta P_n - \sin | \theta | | \hat{P} |) |
\phi_0 |_{\delta}^2 \\
= & \tfrac{1}{2} \int_M \langle (\mu - J \cdot e^0 \cdot) \phi, \phi
\rangle \\
& + \tfrac{1}{2} \int_{\partial M} (H \pm \cos \theta
\ensuremath{\operatorname{tr}}_{\partial M} p) | \phi |^2 \pm \sin \theta
\langle \sqrt{- 1} p_{n \gamma} e^{\gamma} \cdot e^n \cdot e^0 \cdot \phi,
\phi \rangle . \label{mass inequality}
\end{align}
  By the dominant energy conditions \eqref{interior dec} and \eqref{tilt dec},
  we obtain that $E \pm \cos \theta P_n \geqslant \sin | \theta |  | \hat{P}
  |$.
\end{proof}

\begin{lemma}
  There exists a choice of $\phi_0$ and the sign of $\theta$ such that
  \begin{equation}
    \mp \sin \theta P_{\gamma} \langle \sqrt{- 1} \mathrm{d} x^{\gamma} \cdot
    \mathrm{d} x^n \cdot \mathrm{d} x^0 \cdot \phi_0, \phi_0 \rangle_{\delta}
    = \sin | \theta |  | \hat{P} |  | \phi_0 |_{\delta} .
  \end{equation}
\end{lemma}

\begin{proof}
  We work with the flat metric in this proof. Let
  \[ A = P_{\alpha} \sqrt{- 1} \mathrm{d} x^{\alpha} \cdot \mathrm{d} x^n
     \cdot \mathrm{d} x^0, \]
  we know from the last item of Lemma \ref{Q algebraic} that $\bar{Q}$
  commutes with $A$. So they have the same eigen spinors. Then for some
  $\phi_0 \neq 0$ such that $\bar{Q} \phi_0 = \pm \phi_0$, so $A \phi_0 =
  \lambda \phi_0$ for some $\lambda \in \mathbb{C}$. It is easy to check that
  $A$ is Hermitian and so $\lambda \in \mathbb{R}$. We calculate $\lambda$ by
  the following:
\begin{align}
\lambda^2 | \phi_0 |^2 & = P_{\alpha} P_{\beta} \langle \sqrt{- 1}
\mathrm{d} x^{\alpha} \cdot \mathrm{d} x^n \cdot \mathrm{d} x^0 \cdot
\phi_0, \sqrt{- 1} \mathrm{d} x^{\beta} \cdot \mathrm{d} x^n \cdot
\mathrm{d} x^0 \cdot \phi_0 \rangle_{\delta} \\
& = P_{\alpha} P_{\beta} \langle \mathrm{d} x^{\alpha} \cdot \mathrm{d}
x^n \cdot \mathrm{d} x^0 \cdot \phi_0, \mathrm{d} x^{\beta} \cdot
\mathrm{d} x^n \cdot \mathrm{d} x^0 \cdot \phi_0 \rangle_{\delta}
\\
& = P_{\alpha} P_{\beta} \langle \mathrm{d} x^n \cdot \mathrm{d}
x^{\alpha} \cdot \mathrm{d} x^0 \cdot \phi_0, \mathrm{d} x^n \cdot
\mathrm{d} x^{\beta} \cdot \mathrm{d} x^0 \cdot \phi_0 \rangle_{\delta}
\\
& = P_{\alpha} P_{\beta} \langle \mathrm{d} x^{\alpha} \cdot \mathrm{d}
x^0 \cdot \phi_0, \mathrm{d} x^{\beta} \cdot \mathrm{d} x^0 \cdot \phi_0
\rangle_{\delta} \\
& = P_{\alpha} P_{\beta} \langle \mathrm{d} x^0 \cdot \mathrm{d}
x^{\alpha} \cdot \phi_0, \mathrm{d} x^0 \cdot \mathrm{d} x^{\beta} \cdot
\phi_0 \rangle_{\delta} \\
& = P_{\alpha} P_{\beta} \langle \mathrm{d} x^{\alpha} \cdot \phi_0,
\mathrm{d} x^{\beta} \cdot \phi_0 \rangle_{\delta} \\
& = - P_{\alpha} P_{\beta} \langle \mathrm{d} x^{\alpha} \cdot \mathrm{d}
x^{\beta} \cdot \phi_0, \phi_0 \rangle_{\delta} \\
& = - P_{\alpha} P_{\beta} \langle \tfrac{1}{2} (\mathrm{d} x^{\beta}
\cdot \mathrm{d} x^{\alpha} \cdot + \mathrm{d} x^{\alpha} \cdot \mathrm{d}
x^{\beta} \cdot) \phi_0, \phi_0 \rangle_{\delta} \\
& = | \hat{P} |^2 | \phi_0 |_{\delta}^2 .
\end{align}
  Hence $\lambda = \pm | \hat{P} |$. Using this choice of $\phi_0$, from
  \eqref{other pi}, we get
\begin{align}
& \langle P_{\gamma} \mathrm{d} x^{\gamma} \cdot \mathrm{d} x^n \cdot
\mathrm{d} x^0 \cdot \phi_0, \phi_0 \rangle_{\delta} \\
= & \mp \sin \theta \lambda | \phi_0 |_{\delta}^2 \\
= & (\mp 1) (\pm 1) | \hat{P} | \sin \theta | \phi_0 |^2_{\delta} .
\end{align}
  If $| \hat{P} |$ does not vanish, we can always make $\mp \sin \theta
  \lambda$ positive regardless of the sign in $\bar{Q} \phi = \pm \phi$ by
  making a free choice of the sign of $\theta$, fixing such $\theta$, we have
  \[ \langle P_{\gamma} \mathrm{d} x^{\gamma} \cdot \mathrm{d} x^n \cdot
     \mathrm{d} x^0 \cdot \phi_0, \phi_0 \rangle_{\delta} = \sin | \theta |  |
     \hat{P} |  | \phi_0 |_{\delta}^2 . \]
  And the proof is done.
\end{proof}

\subsection{Some consequences of vanishing mass}\label{sec vanishing mass}

We are not able to prove a rigidity statement to Theorem \ref{tilt pmt}.
However, from the equality in Theorem \ref{tilt pmt}, we do have some simple
consequences. We assume that $\theta \neq 0$ in this subsection, however, the
results and the proof work through for the case $\theta = 0$ with only minor
changes.

By \eqref{interior dec}, \eqref{tilt dec}, \eqref{boundary condition} and
\eqref{mass inequality}, there exists a nonzero spinor $\phi$ which satisfies
\begin{align}
\tilde{\nabla} \phi & = 0, \label{parallel} \\
\ensuremath{\operatorname{Re}} \langle (\mu - J \cdot e^0 \cdot) \phi, \phi
\rangle & = 0 \text{ in } M, \label{vanishing constraint} \\
Q \phi & = \pm \phi, \\
\ensuremath{\operatorname{Re}} \langle \tilde{H} \phi, \phi \rangle & = 0
\text{ on } \partial M, \label{vanishing boundary constraint}
\end{align}
where $\tilde{H}$ is a shorthand given by
\begin{equation}
  \tilde{H} = H \pm \cos \theta \ensuremath{\operatorname{tr}}_{\partial M} p
  \mp \sin \theta \sqrt{- 1} p_{\alpha n} e^{\alpha} \cdot e^n \cdot e^0
  \cdot,
\end{equation}
We define
\begin{equation}
  N = \langle \phi, \phi \rangle, X = \sum_j (e^j \cdot \phi, \phi) e_j = X^j
  e_j .
\end{equation}
As the naming of $N$ and $X$ suggests, they are related to the lapse function
and shift vector in Section \ref{sec:invariance}. In the interior of $M$, we
have the following consequence of \eqref{parallel}.

\begin{proposition}
  \label{interior consequence}If \eqref{parallel} holds, then
  \begin{equation}
    L_X g + 2 N p = 0, \text{ } \mathrm{d} (N^2 - |X|^2) = 0
    \label{consequence of parallel}
  \end{equation}
  in $M$. If \eqref{interior dec} and \eqref{vanishing constraint} hold, then
  \begin{equation}
    \mu N + \langle J, X \rangle = 0, \text{ } \mu X^k + N J^k = 0
  \end{equation}
  in $M$.
\end{proposition}

\begin{proof}
  We show by direct calculation. We choose a geodesic normal frame $e_i$ at an
  interior point of $M$, combining with \eqref{parallel}, we have
\begin{align}
e_i (N) & = \tilde{\nabla}_i (e^0 \cdot \phi, \phi) = (\tilde{\nabla}_i
e^0 \cdot \phi, \phi) \\
& = - p_{i j} (e^j \cdot \phi, \phi) = - p_{i j} X^j,
\end{align}
  and
\begin{align}
\nabla_i X & = e_i [(e^j \cdot \phi, \phi) e_j] = (\tilde{\nabla}_i e^j
\cdot \phi, \phi) e_j \\
& = - p_{i j} (e^0 \cdot \phi, \phi) e_j = - p_{i j} N e_j .
\end{align}
  This proves \eqref{consequence of parallel}. From \eqref{vanishing
  constraint}, we obtain
\begin{align}
0 & = \langle (\mu - J \cdot e^0 \cdot) \phi, \phi \rangle \\
& = \mu \langle \phi, \phi \rangle - J_i \langle e^i \cdot e^0 \cdot
\phi, \phi \rangle \\
& = \mu \langle \phi, \phi \rangle + J_i (e^i \cdot \phi, \phi)
\\
& = \mu N + J_i X^i .
\end{align}
  For any $C^2$ spinor $\psi$ and $s \in \mathbb{R}$, by \eqref{interior dec},
  we have
  \begin{equation}
    \ensuremath{\operatorname{Re}} \langle (\mu - J \cdot e^0 \cdot) (\phi + s
    \psi), \phi + s \psi \rangle \geqslant 0.
  \end{equation}
  Subtracting \eqref{vanishing constraint} from the above, we have
\begin{align}
& s\ensuremath{\operatorname{Re}} \langle (\mu - J \cdot e^0 \cdot) \phi,
\psi \rangle + s\ensuremath{\operatorname{Re}} \langle (\mu - J \cdot e^0
\cdot \psi), \phi \rangle \\
+ & s^2 \ensuremath{\operatorname{Re}} \langle (\mu - J \cdot e^0 \cdot)
\psi, \psi \rangle \geqslant 0,
\end{align}
  which forces
  \begin{equation}
    \ensuremath{\operatorname{Re}} \langle (\mu - J \cdot e^0 \cdot) \phi,
    \psi \rangle +\ensuremath{\operatorname{Re}} \langle (\mu - J \cdot e^0
    \cdot \psi), \phi \rangle = 0.
  \end{equation}
  Because $\mu - J \cdot e^0 \cdot$ is Hermitian with respect to $\langle
  \cdot, \cdot \rangle$, so
  \begin{equation}
    \ensuremath{\operatorname{Re}} \langle (\mu - J \cdot e^0 \cdot) \phi,
    \psi \rangle = 0. \label{interior testing}
  \end{equation}
  Setting $\psi$ to be $e^0 \cdot e^k \cdot \phi$ gives
\begin{align}
0 & =\ensuremath{\operatorname{Re}} \langle (\mu - J \cdot e^0 \cdot)
\phi, e^0 \cdot e^k \cdot \phi \rangle \\
& = \mu \langle \phi, e^0 \cdot e^k \cdot \phi \rangle - J_i
\ensuremath{\operatorname{Re}} \langle e^i \cdot e^0 \cdot \phi, e^0 \cdot
e^k \cdot \phi \rangle \\
& = \mu \langle e^0 \cdot e^k \cdot \phi, \phi \rangle + J_i
\ensuremath{\operatorname{Re}} \langle e^i \cdot \phi, e^k \cdot \phi
\rangle \\
& = \mu (e^k \cdot \phi, \phi) + J_i \delta^{i k} | \phi |^2 \\
& = \mu X^k + N J^k .
\end{align}
  And the lemma is proved.
\end{proof}

Now we study the consequences of \eqref{vanishing mass} at the boundary
$\partial M$. First, we show some basic calculations along $\partial M$.

\begin{lemma}
  We have
\begin{align}
\tilde{\nabla}_{\alpha} (Q \phi) = & \cos \theta [- p_{\alpha \beta}
e^{\beta} \cdot e^n \cdot \phi + h_{\alpha \beta} e^0 \cdot e^{\beta}
\cdot \phi] \\
& + \sqrt{- 1} \sin \theta [- p_{\alpha n} e^0 \cdot \phi + h_{\alpha
\beta} e^{\beta} \cdot \phi] . \label{tangential derivative Q phi}
\end{align}
\end{lemma}

\begin{proof}
  First,
  \begin{equation}
    \tilde{\nabla}_i e^0 = - p_{i j} e^j, \tilde{\nabla}_{\alpha} e^n = -
    p_{\alpha n} e^0 + h_{\alpha \beta} e^{\beta} .
  \end{equation}
  By product rule and that $\tilde{\nabla} \phi = 0$, we have
\begin{align}
& \tilde{\nabla}_{\alpha} (\cos \theta e^0 \cdot e^n \cdot \phi + \sqrt{-
1} \sin \theta e^n \cdot \phi) \\
= & \cos \theta [\tilde{\nabla}_{\alpha} e^0 \cdot e^n \cdot \phi + e^0
\cdot \tilde{\nabla}_{\alpha} e^n \cdot \phi] + \sqrt{- 1} \sin \theta
\tilde{\nabla}_{\alpha} e^n \cdot \phi \\
= & \cos \theta [- p_{\alpha i} e^i \cdot e^n \cdot \phi - e^0 \cdot
p_{\alpha n} e^0 \cdot \phi + e^0 \cdot h_{\alpha \beta} e^{\beta} \cdot
\phi] \\
& + \sqrt{- 1} \sin \theta [- p_{\alpha n} e^0 \cdot \phi + h_{\alpha
\beta} e^{\beta} \cdot \phi] \\
= & \cos \theta [- p_{\alpha \beta} e^{\beta} \cdot e^n \cdot \phi + e^0
\cdot h_{\alpha \beta} e^{\beta} \cdot \phi] + \sqrt{- 1} \sin \theta [-
p_{\alpha n} e^0 \cdot \phi + h_{\alpha \beta} e^{\beta} \cdot \phi] .
\end{align}
  \ 
\end{proof}

\begin{lemma}
  We have
  \begin{equation}
    \langle e^{\alpha} \cdot \phi, \phi \rangle = X^{\alpha} \sqrt{- 1} \cot
    \theta .
  \end{equation}
\end{lemma}

\begin{proof}
  From the relation of $Q$ with $e^{\alpha} \cdot e^0 \cdot$ (see Lemma \ref{Q
  algebraic}), we see that
  \begin{equation}
    X^{\alpha} = (e^{\alpha} \cdot \phi, \phi) = \langle e^0 \cdot e^{\alpha}
    \cdot \phi, \phi \rangle = \pm \sqrt{- 1} \sin \theta \langle e^0 \cdot
    e^{\alpha} \cdot e^n \cdot \phi, \phi \rangle . \label{X alpha}
  \end{equation}
  Similarly,
\begin{align}
& \langle e^{\alpha} \cdot \phi, \phi \rangle \\
= & \pm \cos \theta \langle e^{\alpha} \cdot e^0 \cdot e^n \cdot \phi,
\phi \rangle \\
= & - \cos \theta (\pm \langle e^0 \cdot e^{\alpha} \cdot e^n \cdot \phi,
\phi \rangle) \\
= & - \cos \theta \tfrac{X^{\alpha}}{\sqrt{- 1} \sin \theta} = X^{\alpha}
\sqrt{- 1} \cot \theta .
\end{align}
  \ 
\end{proof}

\begin{proposition}
  \label{boundary consequence}We have
  \begin{equation}
    p_{\alpha n} N \sin^2 \theta = \pm p_{\alpha \beta} X^{\beta} \cos \theta
    + h_{\alpha \beta} X^{\beta} = (H \pm \cos \theta
    \ensuremath{\operatorname{tr}}_{\partial M} p) X_{\alpha} \label{boundary
    consequence i}
  \end{equation}
  and
  \begin{equation}
    p_{\alpha n} X^{\gamma} \mp \cos \theta p_{\alpha \gamma} N - h_{\alpha
    \gamma} N = 0 \label{boundary consequence ii}
  \end{equation}
  along $\partial M$.
\end{proposition}

\begin{proof}
  Taking the product $\langle \tilde{\nabla}_{\alpha} (Q \phi), e^0 \cdot
  \phi \rangle$, we get
\begin{align}
& \langle \tilde{\nabla}_{\alpha} (Q \phi), e^0 \cdot \phi \rangle
\\
= & - \cos \theta p_{\alpha \beta} \langle e^{\beta} \cdot e^n \cdot \phi,
e^0 \cdot \phi \rangle \\
& + \cos \theta h_{\alpha \beta} \langle e^0 \cdot e^{\beta} \cdot \phi,
e^0 \cdot \phi \rangle \\
& - \sqrt{- 1} \sin \theta p_{\alpha n} \langle e^0 \cdot \phi, e^0 \cdot
\phi \rangle \\
& + \sqrt{- 1} \sin \theta h_{\alpha \beta} \langle e^{\beta} \cdot \phi,
e^0 \cdot \phi \rangle \\
= & \mp \cos \theta p_{\alpha \beta} \tfrac{X^{\beta}}{\sqrt{- 1} \sin
\theta} \\
& + \cos \theta h_{\alpha \beta} \sqrt{- 1} X^{\beta} \cot \theta
\\
& - \sqrt{- 1} p_{\alpha n} N \sin \theta \\
& + \sqrt{- 1} \sin \theta h_{\alpha \beta} X^{\beta} .
\end{align}
  Since $\tilde{\nabla} \phi = 0$ and $Q \phi = \pm \phi$ along $\partial M$,
  so $\tilde{\nabla}_{\alpha} (Q \phi) = 0$ along $\partial M$. So we get
  \begin{equation}
    p_{\alpha n} N = \pm p_{\alpha \beta} X^{\beta} \tfrac{\cos \theta}{\sin^2
    \theta} + h_{\alpha \beta} X^{\beta} \tfrac{\cos^2 \theta}{\sin^2 \theta}
    + h_{\alpha \beta} X^{\beta}
  \end{equation}
  which leads to the first identity of \eqref{boundary consequence i}. The
  product $\langle \tilde{\nabla}_{\alpha} (Q \phi), e^{\gamma} \cdot \phi
  \rangle$ leads to
\begin{align}
& \langle \tilde{\nabla}_{\alpha} (Q \phi), e^{\gamma} \cdot \phi \rangle
\\
= & - \cos \theta p_{\alpha \beta} \langle e^{\beta} \cdot e^n \cdot \phi,
e^{\gamma} \cdot \phi \rangle + \cos \theta h_{\alpha \beta} \langle e^0
\cdot e^{\beta} \cdot \phi, e^{\gamma} \cdot \phi \rangle \\
& - \sqrt{- 1} \sin \theta p_{\alpha n} \langle e^0 \cdot \phi,
e^{\gamma} \cdot \phi \rangle + \sqrt{- 1} \sin \theta h_{\alpha \beta}
\langle e^{\beta} \cdot \phi, e^{\gamma} \cdot \phi \rangle \\
= & \mp \sin \theta \cos \theta \sqrt{- 1} p_{\alpha \beta} \langle
e^{\gamma} \cdot e^{\beta} \cdot \phi, \phi \rangle \\
& - \sqrt{- 1} \sin \theta p_{\alpha n} X^{\gamma} \\
& - \sqrt{- 1} \sin \theta h_{\alpha \beta} \langle e^{\gamma} \cdot
e^{\beta} \cdot \phi, \phi \rangle .
\end{align}
  Considering $\tilde{\nabla}_{\alpha} (Q \phi) = 0$, we have
  \begin{equation}
    \sin \theta p_{\alpha n} X^{\gamma} \pm \cos \theta \sin \theta p_{\alpha
    \beta} \langle e^{\gamma} \cdot e^{\beta} \cdot \phi, \phi \rangle +
    h_{\alpha \beta} \sin \theta \langle e^{\gamma} \cdot e^{\beta} \cdot
    \phi, \phi \rangle = 0.
  \end{equation}
  Taking the real part and dividing by $\sin \theta$,
  \begin{equation}
    p_{\alpha n} X^{\gamma} \mp \cos \theta p_{\alpha \gamma} N - h_{\alpha
    \gamma} N = 0.
  \end{equation}
  This is \eqref{boundary consequence ii}.
  
  Similar to the derivation of \eqref{interior testing}, we obtain from
  \eqref{tilt dec} and \eqref{vanishing boundary constraint} that for any
  $C^2$ spinor $\psi$ that $\ensuremath{\operatorname{Re}} \langle \tilde{H}
  \phi, \psi \rangle = 0$. Taking $\psi$ to be $e^{\gamma} \cdot \phi$, we
  obtain
\begin{align}
& \langle \tilde{H} \phi, e^{\gamma} \cdot \phi \rangle \\
= & \langle H \phi \pm \cos \theta
\ensuremath{\operatorname{tr}}_{\partial M} p \phi \pm \sin \theta p_{n
\alpha} \sqrt{- 1} e^{\alpha} \cdot e^0 \cdot e^n \cdot \phi, e^{\gamma}
\cdot \phi \rangle \\
= & - (H \pm \cos \theta \ensuremath{\operatorname{tr}}_{\partial M} p)
\langle e^{\gamma} \cdot \phi, \phi \rangle \mp \sqrt{- 1} \sin \theta
p_{\alpha n} \langle e^{\gamma} \cdot e^{\alpha} \cdot e^0 \cdot e^n \cdot
\phi, \phi \rangle \\
= & - (H \pm \cos \theta \ensuremath{\operatorname{tr}}_{\partial M} p)
X^{\gamma} \sqrt{- 1} \cot \theta \mp \sqrt{- 1} \sin \theta p_{\alpha n}
(\pm 1) \cos \theta \langle e^{\gamma} \cdot e^{\alpha} \cdot \phi, \phi
\rangle \\
= & - (H \pm \cos \theta \ensuremath{\operatorname{tr}}_{\partial M} p)
X^{\gamma} \sqrt{- 1} \cot \theta - \sqrt{- 1} \sin \theta p_{\alpha n}
\cos \theta \langle e^{\gamma} \cdot e^{\alpha} \cdot \phi, \phi \rangle .
\end{align}
  So taking the imaginary part of the above, we arrive
  \begin{equation}
    (H \pm \cos \theta \ensuremath{\operatorname{tr}}_{\partial M} p)
    X^{\gamma} + \sin^2 \theta p_{\alpha n} \langle e^{\gamma} \cdot
    e^{\alpha} \cdot \phi, \phi \rangle = 0.
  \end{equation}
  Taking the real part of the above leads to the second identity of
  \eqref{boundary consequence ii}.
\end{proof}

\appendix\section{Stable capillary MOTS}\label{capillary MOTS}

Let $(M, g, p)$ be a compact initial data set, $\gamma$ be a function on $M$
with range $(0, 2 \pi)$. Let $\Sigma$ be a hypersurface in $M$ with boundary
$\partial \Sigma \subset \partial M$.

\begin{figure}[h]
  \resizebox{102pt}{84pt}{\includegraphics{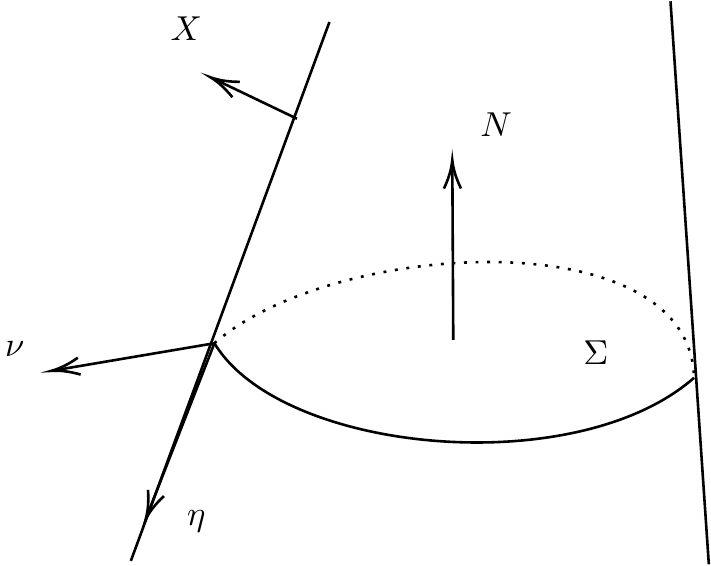}}
  \caption{Naming of various vectors.\label{naming of vectors}}
\end{figure}

We fix some notations used in this Appendix:
\begin{itemize}
  \item $N$, unit normal of $\Sigma$ in $M$; $A$, the second fundamental form
  of $\Sigma$ in $M$; $H$, the mean curvature of $\Sigma$ in $M$;
  
  \item $\eta$, unit normal of $\partial \Sigma$ in $\partial M$;
  
  \item $X$, unit outward normal of $\partial M$ in $M$; $A_{\partial M}$, the
  second fundamental form of $\partial M$ in $M$; $H_{\partial M}$, the mean
  curvature of $\partial M$ in $M$;
  
  \item $\nu$, unit outward normal of $\partial \Sigma$ in $\Sigma$;
  $H_{\partial \Sigma}$, the mean curvature of $\partial \Sigma$ in $\Sigma$.
\end{itemize}
See the Figure \ref{naming of vectors}.

\begin{definition}
  We say that a hypersurface $\Sigma \subset M$ with boundary $\partial \Sigma
  \subset \partial M$ is a marginally trapped surface (in short, MOTS) if
  \begin{equation}
    H +\ensuremath{\operatorname{tr}}_M p = 0. \label{mots}
  \end{equation}
  It is called capillary if
  \[ \langle X, N \rangle = \cos \gamma \label{contact angle} \text{ along }
     \partial \Sigma . \]
\end{definition}

\begin{definition}
  \label{def:stable capillary MOTS}We say that $\Sigma$ is stable capillary
  MOTS if it satisfies both \eqref{mots} and \eqref{contact angle}, and
  furthermore there exists a vector field $Y$ tangent to $\partial M$ along
  $\partial \Sigma$ such that $\langle Y, N \rangle > 0$,
  \begin{equation}
    \nabla_Y (H +\ensuremath{\operatorname{tr}}_M p) \geqslant 0 \text{ in }
    \Sigma
  \end{equation}
  and
  \begin{equation}
    \nabla_Y (\langle X, N \rangle - \cos \gamma) = 0 \text{ along } \partial
    \Sigma . \label{eq:boundary stability condition}
  \end{equation}
\end{definition}

We generalize \eqref{tilt dec} slightly by allowing varying angles as follows.

\begin{definition}
  We say that $(M, g, p)$ satisfies the generalized tilted boundary dominant
  energy condition if the boundary satisfies the inequality
  \begin{equation}
    H_{\partial M} + \cos \gamma \ensuremath{\operatorname{tr}}_{\partial M} p
    \geqslant \sin \gamma |p (X, \cdot) - \tfrac{1}{\sin \gamma} \langle
    \nabla^{\partial M} \gamma, \cdot \rangle | . \label{capillary dec}
  \end{equation}
\end{definition}

We follow {\cite{galloway-generalization-2006}} to obtain the following. It is
also a generalization of {\cite{mendes-rigidity-2022}} where free boundary
MOTS were treated.

\begin{theorem}
  If $M$ satisfies the dominant energy condition \eqref{interior dec} and
  \eqref{capillary dec}, then a stable capillary MOTS is of positive Yamabe
  type unless $\Sigma$ is Ricci flat with totally geodesic boundary.
\end{theorem}

\begin{proof}
  Let $Y$ be a vector field given in Definition \ref{def:stable capillary
  MOTS} and $\phi = \langle Y, N \rangle$. By definition, $\phi$ is positive.
  Using the variation formula for $H +\ensuremath{\operatorname{tr}}_{\Sigma}
  p$ ({\cite{galloway-generalization-2006}}, {\cite{andersson-local-2005}} or
  {\cite[Proposition 2]{eichmair-spacetime-2016}}), we have that
  \[ L \phi := - \Delta \phi + 2 \langle W, \nabla \phi \rangle + (Q
     +\ensuremath{\operatorname{div}}W - |W|^2) \phi \geqslant 0. \]
  where $Q = \tfrac{1}{2} R_{\Sigma} - (\mu + J (N)) - \tfrac{1}{2} | \chi
  |^2$, $\chi = A + p$, $W$ is the vector field tangential to $\Sigma$ dual to
  the 1-form $p (N, \cdot)$ and $R_{\Sigma}$ is the scalar curvature of
  $\Sigma$. Let $Z = W - \nabla \log \phi$, as in
  {\cite{galloway-generalization-2006}}, we get
  \[ - Q \psi^2 \leqslant \ensuremath{\operatorname{div}} (\psi^2 Z) + |
     \nabla \psi |^2 \]
  for any $\psi \in C^{\infty} (\Sigma)$. We integrate the above over $\Sigma$
  and using the divergence theorem, we have that
  \begin{equation}
    \int_{\partial \Sigma} \psi^2 \langle Z, \nu \rangle + \int_{\Sigma} |
    \nabla \psi |^2 + Q \psi^2 \geqslant 0. \label{symmetrized stable}
  \end{equation}
  Now we compute
\begin{align}
\langle Z, \nu \rangle = & \langle W, \nu \rangle - \langle \nabla \log
\phi, \nu \rangle \\
= & p (N, \nu) - (\tfrac{1}{\sin \gamma} A_{\partial M} (\eta, \eta) -
\cot \gamma A (\nu, \nu) + \tfrac{1}{\sin^2 \gamma} \partial_{\eta} \cos
\gamma) \\
= & p (N, \nu) - (- H \cot \gamma + \tfrac{H_{\partial M}}{\sin \gamma} -
H_{\partial \Sigma}) + \tfrac{1}{\sin \gamma} \partial_{\eta} \gamma
\\
= & H_{\partial \Sigma} - \tfrac{1}{\sin \gamma} (H_{\partial M} - H \cos
\gamma - \sin \gamma p (N, \nu)) + \tfrac{1}{\sin \gamma} \partial_{\eta}
\gamma \\
= & H_{\partial \Sigma} - \tfrac{1}{\sin \gamma} (H_{\partial M}
+\ensuremath{\operatorname{tr}}_{\Sigma} p \cos \gamma - \sin \gamma p (N,
\nu) - \partial_{\eta} \gamma) .
\end{align}
  In the above, the second identity follows from computing \eqref{eq:boundary
  stability condition} explicitly whose details could be found in the Appendix
  of {\cite{ros-stability-1997}}. The third identity could be found in
  {\cite[Lemma 3.1]{ros-stability-1997}}. We claim that
  \[ \ensuremath{\operatorname{tr}}_{\Sigma} p \cos \gamma - \sin \gamma p (N,
     \nu) = \cos \gamma \ensuremath{\operatorname{tr}}_{\partial M} p + \sin
     \gamma p (X, \eta) . \]
  Let $\{e_i \}_{1 \leqslant i \leqslant n}$ be an orthonormal frame on $T
  \Sigma$ such that $\nu = e_n$. Then
\begin{align}
& \ensuremath{\operatorname{tr}}_{\Sigma} p \cos \gamma - \sin \gamma p
(N, \nu) \\
= & \cos \gamma \sum_{i \neq n} p_{i i} + \cos \gamma p_{\nu \nu} - \sin
\gamma p_{N \nu} \\
= & \cos \gamma \sum_{i \neq n} p_{i i} + p (\cos \gamma \nu - \sin \gamma
N, \nu) \\
= & \cos \gamma \sum_{i \neq n} p_{i i} + p (\eta, \nu) \\
= & \left( \sum_{i \neq n} \cos \gamma p_{i i} + \cos \gamma p_{\eta \eta}
\right) + p (\eta, \nu - \cos \gamma \eta) \\
= & \cos \gamma \ensuremath{\operatorname{tr}}_{\partial M} p + \sin
\gamma p (\eta, X)
\end{align}
  From the inequality \eqref{symmetrized stable}, we have
\begin{align}
& \int_{\Sigma} R_{\Sigma} \psi^2 + 2 \int_{\partial \Sigma} H_{\partial
\Sigma} \psi^2 \\
\geqslant & \int_{\Sigma} 2 (\mu - J (N)) \psi^2 + \int_{\partial \Sigma}
\tfrac{1}{\sin \gamma} [H_{\partial M} + \cos \gamma
\ensuremath{\operatorname{tr}}_{\partial M} p + \sin \gamma p (\eta, X) -
\partial_{\eta} \gamma] .
\end{align}
  By the dominant energy conditions \eqref{interior dec} and \eqref{capillary
  dec}, we see
  \[ \int_{\Sigma} R_{\Sigma} \psi^2 + 2 \int_{\partial \Sigma} H_{\partial
     \Sigma} \psi^2 \geqslant 0 \]
  confirming that $\Sigma$ is of nonnegative Yamabe type. The rest concerns an
  eigenvalue problem on $\Sigma$, it is no different from the free boundary
  case, see {\cite{mendes-rigidity-2022}}.
\end{proof}

\begin{remark}
  The quantity $H +\ensuremath{\operatorname{tr}}_{\Sigma} p$ in \eqref{mots}
  is called \text{{\itshape{outer null expansion}}}. Further results such as
  constructing constant null expansion foliation and related rigidity analysis
  are not difficult, the readers may again consult
  {\cite{mendes-rigidity-2022}}.
\end{remark}

\

\end{document}